\documentclass[12pt ]{amsart}
\usepackage[pdftex]{graphicx}
\usepackage{amsmath,amssymb,amsthm,mathrsfs,ascmac,comment,lineno}
\usepackage[top=30truemm,bottom=30truemm,left=30truemm,right=30truemm]{geometry}
\usepackage[all]{xy}
\usepackage[capitalize]{cleveref}
\usepackage{float}
\usepackage{caption}

\newcommand{\ctext}[1]{\raise0.2ex\hbox{\textcircled{\scriptsize{#1}}}}

\newcommand{\Q}{\mathbb{Q}}

\newcommand{\C}{\mathbb{C}}
\newcommand{\B}{\mathbb{B}}
\newcommand{\Z}{\mathbb{Z}}
\newcommand{\F}{\mathbb{F}}
\newcommand{\e}{\epsilon}

\newcommand{\s}{\sigma}

\newcommand{\Gal}{\textup{Gal}}
\newcommand{\Cl}{\textup{Cl}}

\theoremstyle{plain}
\newtheorem{defi}{Definition}[section]
\newtheorem{theo}[defi]{Theorem}
\newtheorem{prop}[defi]{Proposition}
\newtheorem{lem}[defi]{Lemma}
\newtheorem{cor}[defi]{Corollary}

\newtheorem{prob}[defi]{Problem}
\newtheorem{conj}[defi]{Conjecture}

\title[Weber's class number problem and its variants]{Weber's class number problem and its variants}
\dedicatory{Dedicated to Professor Toshie Takata}

\makeatletter
\@namedef{subjclassname@2020}{
\textup{2020} Mathematics Subject Classification}
\makeatother 

\author{Hyuga Yoshizaki}
\email{yoshizaki.hyuga@gmail.com}
\address{Department of Mathematics, Graduate school of science and Technology, Tokyo University of Science, 2641, Yamazaki, Noda-shi, 278-8510, Chiba, Japan}

\subjclass[2020]{Primary
11R29, 
11R58, 
57M12; 
Secondary
11E95, 
11R23, 
11R27, 
14H52, 
57K14. 
}

\keywords{Weber's class number problem, arithmetic topology, Iwasawa theory, Pell's equation, units, elliptic curves, knots.}

\date{}

\begin{document}

\begin{abstract}
We survey Weber's class number problem and its variants in the spirit of arithmetic topology;
we recollect some history, present a relation to certain units and generalized Pell's equation,
and overview {a} study of the $p$-adic limits of class numbers in $\Z_p$-towers together with numerical investigation for knots and elliptic curves.
\end{abstract}
\maketitle


\tableofcontents

\section{Introduction --Weber's class number problem--}

Let $k$ be a number field, that is, a finite extension over the rational number field $\Q$.
{Let $\mathcal{O}_k$ denote the ring of integers of $k$, that is, the set of roots of integer coefficients monic polynomials in $k$.
The ideal class group of $k$ is a quotient of the group of fractional ideals of $\mathcal{O}_k$ by the group of principal ideals.}
The class group is finite, and its size $h(k)$ is called the class number of $k$.
The class number is an important invariant of a number field and is a central object of interest in number theory.
For example, Gauss conjectured that there are infinitely many real quadratic number fields with class number $1$, but that is still open.
It is also unsolved whether there are infinitely many number fields with class number $1$ or not.

Weber's class number problem (conjecture) asks whether {a certain infinitely many number fields} have the class numbers $1$.
Let $p$ be a prime number.
Let $\Z_p$ denote the ring of $p$-adic integers.
A $\Z_p$-extension over a number field is an infinite Galois extension whose Galois group is isomorphic to the additive group $\Z_p$ as topological groups.
The statement of Weber's problem is as follows.
  \begin{prob}
    Determine the class number of each intermediate field in the $\Z_p$-extension over $\Q$.
    Is the class number of any intermediate field $1$?
  \end{prob}
Many researchers expect that the answer to this problem is affirmative.

This article is a survey based on several articles containing \cite{Yoshizaki} and \cite{Ueki-Yoshizaki}.
In the former half of this article, we survey number theory related to Weber's problem. 
We recollect two approaches and a recent trial using generalized Pell's equation and certain units. 
In addition, we state a substantial result, that is, the $p$-adic convergence of class numbers. 

{In the latter half of the article, we discuss analogous studies (problems) in the contexts of function fields and compact $3$-manifolds.
Besides Weber's problem is unsolved in number fields, its analogues are solved in function fields and compact  $3$-manifolds.
Analogues of the $p$-adic convergence theorem also hold.
Furthermore, we can compute the limits concretely for some cases.
We position the studies of the $p$-adic limit as variants of Weber's problem, and investigated numerically the $p$-adic limit for knots and elliptic curves.}
For this purpose, for $f(t) \in \Z[t]$, we invoke results on the cyclic resultants $\text{Res}(t^n-1, f(t)) =\prod_{\zeta^n=1} f(\zeta)$ and their $p$-adic limits.

The contents of this article is as follows.
In Section 2, we review the history of Weber's problem and recollect two approaches and present a recent trial.
In particular, we introduce the $p$-adic limit of class numbers which is a main target of the later half of this article.
In Section 3, we overview the {analogous results to number fields for} function fields and compact $3$-manifolds.
In Section 4, we give a numerical investigation of the $p$-adic limits for knots and elliptic curves.
Some concrete examples are different from those in \cite{Ueki-Yoshizaki}.

\section{Number fields}
In this section, we focus on {the most well-studied case} $p=2$.
{We present} two approaches and a recent trial to Weber's problem.

\subsection{History of Weber's problem}
For each positive integer $n$, let $\B_n$ denote the $n$-th layer {of the $\Z_2$-extension over $\Q$}, that is, the unique intermediate field of degree $2^n$.
In 1886, H. Weber \cite[Theorem C]{HWeber} showed that $h(\B_n)$ are odd for all positive integers $n$.
He also showed $h(\B_n)=1$ for $n=1$, $2$, and $3$ by hand calculations.

After Weber's study,
the development of computers allowed researchers to {determine} the class numbers of $\B_n$ as bellow.
\begin{itemize}
  \item $h(\B_4)=1$ (Bauer \cite[Ergebnis]{Bauer}, Masley \cite[Theorem 3.2]{Masley})
  \item $h(\B_5)=1$ (Linden \cite[Theorem 1]{Linden})
  \item $h(\B_6)=1$ (Miller \cite[Theorem 2.1]{Miller})
  \item $h(\B_7)=1$ under Generalized Riemann Hypothesis (Miller \cite[Theorem 2.2]{Miller})
\end{itemize}
Except for Bauer, they first gave an upper bound of the class number $h(\B_n)$ by using {\it the root discriminant} for each $n$ and
then showed that $h(\B_n)$ is not divisible by prime numbers below that upper bound.
The techniques and their difficulties are briefly summarized in the introduction of \cite{Miller}.
{We count them as the first approach.}

The second approach {we present here is to determinate} prime numbers that do not divide $h(\B_n)$ for all $n$.
Fukuda beautifully summarized this approach in his book \cite[Chapter 14]{Fukuda-Kaisetsu}, which is written in Japanese.
See also \cite{FKM} for similar contents shortly written in English.
The explanation {below} is based on his book.

Horie found that there is a strong relation between a prime number dividing $h(\B_n)$ and a certain unit in $\B_n$.
We need some preparations to explain Horie's work.
We write $\zeta_{2^n}=\exp(2\pi\sqrt{-1}/2^n)$ for each $n\geq1$.
Take a generator $\sigma$ of $\Gal(\Q(\zeta_{2^{n+2}})/\Q(\zeta_{2^2}))\cong \Z/2^n\Z$.
For $\alpha=\sum_{i=0}^{2^{n-1}-1}a_i\zeta_{2^n}^i\in \Z[\zeta_{2^n}]$ ($a_i\in \Z$), we define
\[
  \alpha_{\sigma}=\sum_{i=0}^{2^{n-1}-1}a_i\sigma^i \in \Z[\Gal(\Q(\zeta_{2^{n+2}})/\Q(\zeta_{2^2}))].
\]
{For a prime number $l$,} let $F_l$ be the decomposition field of $l$ in the extension $\Q(\zeta_{2^n})/\Q$ and fix an intermediate field $F$ of $\Q(\zeta_{2^n})/F_l$.
Define
\[
 \eta_n=\tan\frac{\pi}{2^{n+2}} 
\]
for each $n\geq1$. The following Horie's lemma plays a key role in this approach.
\begin{lem}[cf.~{\cite[Lemma 2]{Horie05}}]\label{Horie}
 A prime number $l$ divides $h(\B_n)/h(\B_{n-1})$ if and only if
 there is a prime ideal $\mathfrak{L}$ of $F$ dividing $l$ such that $\sqrt[l]{\eta_n^{\alpha_\sigma}}\in \B_n$ for all $\alpha\in l\mathfrak{L}^{-1}$.
\end{lem}
By showing that ``$l\mid h(\B_n)/h(\B_{n-1})$ implies $\sqrt[l]{\eta_n^{\alpha_\sigma}}\in \B_n$'' and estimating a particular size of $\eta_n$, Horie {deduced} a contradiction.
He obtained the following.
\begin{theo}[{\cite[Theorem 3]{Horie07}}]
  If a prime number $l$ satisfies $l\not\equiv \pm1 \pmod 8$, then we have $l\nmid h(\B_n)$ for all $n\geq1$.
\end{theo}
Many researchers have improved Horie's result.

On the other hand, Fukuda and Komatsu \cite[Theorem 1.2]{Fukuda-Komatsu1} showed that for each prime number $l$ there exists an integer $m_l$ such that
$l\nmid h(\B_{m_l})$ {implies} $l\nmid h(\B_n)$ for all $n\ge1$.
They also improved their result in \cite[Theorem 1.1]{Fukuda-Komatsu3}.
By computing such $m_l$ and showing $l\nmid h(\B_{m_l})$ for small $l$, they obtained the following.
\begin{theo}[{\cite[Corollary 1.2]{Fukuda-Komatsu3}}]\label{smallprime}
  If a prime number $l$ satisfies $l<10^9$, then we have $l\nmid h(\B_n)$ for all $n\geq1$.
\end{theo}

Consolidating the recent results, we can summarize as following.
\begin{theo}[{\cite[Corollary B]{Morisawa-Okazaki}}]
  If a prime number $l$ satisfies $l\not\equiv \pm1 \pmod {64}$, then we have $l\nmid h(\B_n)$ for all $n\geq1$.
\end{theo}

\subsection{Units}
{As we have seen above}, it is indispensable to study the group of units of $\B_n$ for Weber's problem.
In this section, we focus on the group of units of $\B_n$.
Let $E_n$ denote the group of units of $\B_n$.
The group of cyclotomic units is defined by
\[
  C_n=\langle -1, \zeta_{2^{n+2}}^{-1}\frac{1-\zeta_{2^{n+2}}^3}{1-\zeta_{2^{n+2}}} \rangle_{\Z[\Gal(\B_n/\Q)]},
\]
where $\langle\cdot\rangle_{\Z[\Gal(\B_n/\Q)]}$ denote {the subgroup generated as a $\Gal(\B_n/\Q)$-module (cf.~\cite[Lemma 8.1, Proposition 8.11]{Washington}).}
Then we have 
\begin{equation}\label{hEC}
  h(\B_n)=\#(E_n/C_n)
\end{equation}
for all $n\geq1$ (see \cite[Theorem 8.2]{Washington}).
We remark that the proof of \cref{Horie} is based on this result.
The equality (\ref{hEC}) implies that it {suffices} to determine fundamental units of $E_n$ to resolve Weber's problem.

In fact, it suffices to determine generators of a certain subgroup of $E_n$.
We define a map ${\rm N}_{n/n-1}:E_n\to E_{n-1}$ by $\e\mapsto \e\tau(\e)$, where $\tau$ is the generator of $\Gal(\B_n/\B_{n-1})\cong \Z/2\Z$.
Set
\[
  RE_n^+=\ker {\rm N}_{n/n-1}=\{\e\in E_n\mid {\rm N}_{n/n-1}(\e)=1\}.
\]
We note that $RE_n^+$ is a subgroups of $E_n$, and $\Gal(\B_n/\Q)$ acts on it.
Since ${\rm N}_{n/n-1}$ is surjective on $E_n$ and $C_n$, we obtain the exact sequence
\begin{equation}\label{exact}
  1\rightarrow RE_n^+/A_n \rightarrow E_n/C_n \rightarrow E_{n-1}/C_{n-1} \rightarrow 1,
\end{equation}
where {we put} $A_n=RE_n^+\cap C_n$.
Hence we have $h(\B_n)/h(\B_{n-1})=(RE_n^+:A_n)$.

There is an interesting conjecture on $RE_n^+$ relating the second approach.
We define $RE_n={\rm N}_{n/n-1}^{-1}(\{\pm1\})$ and $RE_n^-={\rm N}_{n/n-1}^{-1}(\{-1\})$.
Then we have $RE_n=RE_n^+\sqcup RE_n^-$.
For $\e\in RE_n$, we define ${\rm Tr}_{n}(\e)=\sum_{\s \in \Gal(\B_n/\Q)}\s(\e)$.
In the personal communication between Komatsu, Morisawa and Okazaki, they stated the following conjecture.
\begin{conj}\label{Komatsu-Okazaki}
  For all $n\geq1$, we have
  \[
    \min\{{\rm Tr}_n(\e^2)\mid \e \in RE_n\setminus \{\pm1\}\}\geq 2^n(2^{n+1}-1).
  \]
\end{conj}
Morisawa and Okazaki studied this conjecture and resolved for $RE_n^-$ in \cite[Theorem 6.4]{MO}.
Fukuda and Komatsu used their result for \cite[Theorem 1.3]{Fukuda-Komatsu3}.
In contrast, \cref{Komatsu-Okazaki} has not yet been resolved for $RE_n^+$.
For $RE_n^+$, Morisawa and Okazaki obtained the following.
\begin{prop}[{\cite[Proposition 6.6]{MO}}]
  For all $n\geq2$, we have
  \[
    \min\{{\rm Tr}_n(\e^2)\mid \e \in RE_n^+\setminus \{\pm1\}\}\geq 2^n\cdot 17.
  \]
\end{prop}
Moreover, they extracted informations of $RE_n^+$ by combining two filtrations on $RE_n^+$ (cf.~\cite[Subsection 6.2]{MO}).

Kashio and the author refined \cref{Komatsu-Okazaki} in \cite{Kashio-Yoshizaki}.
\begin{conj}[{refinement of \cref{Komatsu-Okazaki}, {\cite[Conjecture 2.2]{Kashio-Yoshizaki}}}]\label{refinement}
  Let $c_1=2$, $c_n=2\cdot \textup{round}(2^n/5)$ ($n\geq2$) where $\textup{round}(x)$ denotes the nearest integer to $x$.
  Then, for all $n\geq1$, we have
  \[
    \min\{{\rm Tr}_n(\e^2)\mid \e \in RE_n^+\setminus \{\pm1\}\}= 2^n(1+8c_n).
  \]
\end{conj}
For the refined conjecture, we obtained the following.
\begin{theo}
  \cref{refinement} holds for $1\leq n \leq 6$.
\end{theo}
This theorem is proved by a combinatorial argument for $1\leq n \leq 3$.
For the cases $4\leq n \leq 6$, we use $h(\B_n)=1$ and check the conjecture by numerical calculations by using PARI/GP.

\subsection{Generalized Pell's equation}
Recently, the author studied the group $RE_n^+$ from a new approach.
Set $X_n=2\cos(2\pi/2^{n+2})$ for each $n\geq0$. Then we have
\[
  X_1=\sqrt{2},\ \ X_2=\sqrt{2+\sqrt{2}},\ \ X_3=\sqrt{2+\sqrt{2+\sqrt{2}}},....
\]
Since $\B_n$ is the maximal real subfield of $\Q(\zeta_{2^{n+2}})$, we see that the ring of integers of $\B_n$ is $\Z[X_n]$ (cf.~\cite[Propositon 2.16]{Washington}).
For each $\e \in RE_n^+$, there exist unique algebraic integers $a,b \in \Z[X_{n-1}]$ such that $\e=a+X_nb$.
Then we have ${\rm N}_{n/n-1}(\e)=a^2-X_n^2b^2$.
{Thus, the determination of $RE_n^+$ reduces to finding} the all $\Z[X_{n-1}]$-solutions of the equation
\begin{equation}\label{pell}
  x^2-X_n^2y^2=1.
\end{equation}
By sending $a+X_nb$ to $(a,b)$, we have a bijection between $RE_n^+$ and the solutions of \cref{pell}
and identify these sets by this bijection.

In the case of $n=1$, $x^2-2y^2=1$, this equation is called Pell's equation.
In general, for a non-square positive integer $m$, the equation $x^2-my^2=1$ is called Pell's equation.
By Dirichlet's unit theorem, there is a generator of $RE_1^+$, and we call the corresponding solution a fundamental solution.
There is a classical algorithm to find a fundamental solution of Pell's equation by using the regular continued fraction expansion of $\sqrt{m}$.
We note that the algorithm is based on the theory of approximation, which is called the best approximation theorem.

Our strategy for solving \cref{pell} is to imitate the aforementioned classical algorithm.
{By generalizing the classical continued fraction expansion algorithm, we obtain a new continued fraction expansion of $X_n$ over $\Z[X_{n-1}]$.}
\begin{theo}[{\cite[Theorem 3.4]{Yoshizaki}}]
  For each $n\geq1$, we have
  \[
  X_n=[1,\overline{2(1+X_{n-1})^{-1},2}]=1+\cfrac{1}{2(1+X_{n-1})^{-1}+\cfrac{1}{2+\ddots}}.
  \]
\end{theo}
From this continued fraction, by imitating the classical algorithm, we obtain the explicit solution for each $n\geq1$.
{Let $\e_n$ denote the corresponding unit in $RE_n^+$.
Then we have}
\[
  \e_n=\frac{X_n+1}{X_n-1}.
\]
Here we state
\begin{conj}\label{conj}
  $RE_n^+=\langle -1,\e_n \rangle_{\Gal(\B_n/\Q)}$  for all $n\geq1$.
\end{conj}

On the other hand, we have that $A_n(=RE_n^+\cap C_n)=\langle -1,\e_n \rangle_{\Gal(\B_n/\Q)}$ (see \cite[Theorem 4.1]{Yoshizaki} for details).
Thus we see that 
Weber's class number problem is in affirmative if and only if \cref{conj} holds.
{Further discussion for \cref{conj} is in \cite[Section 6]{Yoshizaki}.}

\subsection{$p$-adic limits}
We shall remark some applications of the explicit unit $\e_n\in A_n$.
By using $\e_n$ and counting the orbital decomposition of each prime part of $RE_n^+/A_n$ with respect to the action of $\Gal(\B_n/\Q)$, we obtain
\begin{prop}[cf.~{\cite[Theorem 5.3]{Yoshizaki}}]\label{Ylimit}
  $\frac{h(\B_n)}{h(\B_{n-1})}\equiv 1 \pmod{2^n}$ for each $n\ge1$.
\end{prop}
Since $\Z_2$ is a complete metric space in $2$-adic topology, we see that the sequence of the class numbers $\{h(\B_n)\}_n$ converges in $\Z_2$.
The following theorem is a generalization of \cref{Ylimit}.

\begin{theo}[{\cite[Corollary 2]{Kisilevsky}, \cite[Theorem 2.1]{Ueki-Yoshizaki}}] \label{nconv} 
  Let $k_{p^\infty}$ be a $\Z_p$-extension over a global field $k$,
  and $k_{p^n}$ be the $n$-th layer.
  Then the sizes of the class groups $\Cl(k_{p^n})$, 
  those of the non-$p$-subgroups $\Cl(k_{p^n})_{\text{non-}p}$, 
  and those of the $l$-torsion subgroups $\Cl(k_{p^n})_{(l)}$ for each prime number $l$ 
  converge in $\Z_p$. 
\end{theo}
{This assertion is proved by using ideal class groups, instead of unit groups.}
We remark that we obtained this result independently from Kisilevsky, and the proofs are different.
An essential fact for our proof is the existence of a surjective system of the $l$-parts of the class groups $\Cl(k_{p^n})_{(l)}\to \Cl(k_{p^{n-1}})_{(l)}$ for each prime number $l\neq p$.
We take the relative norm map of number fields, which is defined by $N_{k_{p^n}/k_{p^{n-1}}}(x)=\prod_{\tau \in \Gal(k_{p^n}/k_{p^{n-1}})}\tau(x)$.

{Though the $p$-adic limits of class numbers} exist for all $\Z_p$-extensions over number fields, we have no idea how to compute them.
A strategy to study these invariants may be to study the p-adic limits in analogous contexts, say, that of function fields and $3$-manifolds. 

\section{Analogues}
In this section, based on \cite{Ueki-Yoshizaki},
we discuss analogous results to number theory in the contexts of function fields and $3$-manifolds.

\subsection{Function fields}\label{Inff}
Let $\F$ be a finite field and $t$ be a transcendental element of $\F$.
We call a finite extension $k$ of $\F(t)$ a function field. Set $\F_k=k\cap \bar{\F}$, where $\bar{\F}$ is an algebraic closure of $\F$.
An {extension $k'/k$ over a function field is} said to be constant if $k'=\F_{k'}k$ and geometric if $\F_{k'}=\F_k$.

{There is a deep analogy} between number fields and function fields (cf.~\cite{Iwasawa} in Japanese),
and the theories of number fields and that of function fields have been mutually developed (cf.~\cite{Rosen}).
For example, the analogue of the ideal class group in a function field $k$ is the degree $0$ divisor class group $\Cl^0(k)$, and they are known to be finite (see \cite[Lemma 5.6]{Rosen}).
We call the order of $\Cl^0(k)$ the class number of $k$ and denote it by $h(k)$.

Arguments for function fields often become easier than those for number fields.
The following theorem {states} an answer to a class number problem in function fields.

\begin{theo}[{\cite[Theorem 1.1]{Mercuri-Stirpe}}, {\cite{Shen-Shi}}]
  Except for rational function fields, there are exactly $8$ function fields of class number $1$ up to isomorphism.
\end{theo}
From this result, we find that Weber's conjecture in function fields does not hold except for constant $\Z_p$-extensions over rational function fields.

Unlike in the case of number fields, we can compute the $p$-adic limits for constant $\Z_p$-extensions over function fields.
We must remark that Kisilevsky also obtained the same result and compute the $p$-adic limit.

We position a study of the $p$-adic limits as a variant of Weber problem.
For example, {it is a natural question to ask when the} $p$-adic limits are $1$.
We will answer this question for function fields of genus $1$ in \cref{Examples}.

As a preparation for computation of $p$-adic limits,
we explain the class number formula of function fields in the following.
Let $k$ be a function field with a constant field $\F_k$ of $q$ elements.
The congruent zeta function of $k$ is defined by
\begin{equation}
  \zeta_k(s)=\sum_{A:\text{effective divisor}}\frac{1}{(q^{\deg A})^s}.
\end{equation}
It is known that this function absolutely converges to a holomorphic function on ${\rm Re}(s)>1$.
Moreover, we have the following.
\begin{prop}[{Hasse--Weil, cf.~\cite{Weil}, \cite[Theorem 5.9]{Rosen}}] There exists $L_k(t)\in \Z[t]$ of degree $2g(k)$ satisfying  
\[\zeta_k(s)=\dfrac{L_k(q^{-s})}{(1-q^{-s})(1-q^{1-s})}\]
on ${\rm Re}(s)>1$. 
\end{prop} 
The right-hand side is an analytic continuation of $\zeta_K(s)$ to $\mathbb{C}$ as a meromorphic function.
The polynomial $L_k(t)$ is called the $L$-polynomial of $k$.
The $L$-polynomial satisfies a reciprocal property $L_k(t)=q^{g(k)}t^{2g(k)}L_k(1/(qt))$
and is decomposed by $L_k(t)=\prod_{i=1}^{2g(k)}(1-\alpha_i t)$ for some algebraic integers $\alpha_i$ with $|\alpha_i|=\sqrt{q}$.
In addition, the class number $h(k)$ coincides with $L_k(1)$.
If $k_n/k$ is a constant extension of degree $n$, then we have $k_n=\F_{k_n}k$, $g(k_n)=g(k)$, and $L_{k_n}(t)=\prod_{i=1}^{2g(k)}(1-\alpha_i^n t)$.
Therefore,
\[h(k_n)=L_{k_n}(1)=\prod_{i=1}^{2g(k)}(1-\alpha_i^n)={\rm Res}(t^n-1, t^{2g(k)}L_k(1/t)).\]
We note that there is {a natural bijection} between the isomorphism classes of function fields over a finite field $\F$
and the isomorphism classes of nonsingular projective algebraic curves over $\F$.
For any prime number $l$ relatively prime to $q$, 
the Frobenius polynomial $F_k(t)$ of $k$ is defined as 
the characteristic polynomial of the geometric Frobenius action on the $l$-adic \'{e}tale cohomology of the algebraic curve corresponding to $k$ (cf.~\cite{Aubry-Perret}).
Since we have $F_k(t)=t^{2g(k)}L_k(1/t)$, we obtain the following. 

\begin{prop} \label{FoxWeberFF} 
Let $k$ be a function field and $k'/k$ a constant extension of degree $n$. Then 
\[h(k')=|{\rm Res}(t^n-1,F_k(t))|.\]
\end{prop} 
This formula may be seen as an analogue of Fox--Weber's formula (\cref{FoxWeber}) for a constant extension over a function field. 
{Based on this formula, we give some concrete computations for the $p$-adic limits of function fields of elliptic curves in Subsection \ref{ExamplesFF}.}

To the contrary, in the case of geometric $\Z_p$-extensions, the situation is more complicated,
and we have not yet succeeded in computing the $p$-adic limits.

\subsection{Compact $3$-manifolds}\label{Intop}
{The analogy between prime numbers and knots has contributed to the development of each theory and raised issues with each other.}
We focus on the analogy between the Alexander--Fox theory for $\Z$-covers and the Iwasawa theory for $\Z_p$-extensions, which has played an important role since 1960's (cf.~\cite{Mazur}).

Let $X$ be a compact $3$-manifold.
The first homology group of $\Z$-coefficient $H_1(X)$ is an analogue of the ideal class group,
and the size of its torsion part $|H_1(X)_{\textup{tor}}|$ is an analogue of the class number.
A $\Z_p$-cover of $X$ is a compatible system $(X_{p^n}\to X)_n$ of $\Z/p^n\Z$-covers.

In particular, we mainly study the $p$-power-th cyclic coverings on $S^3$ branched along a knot.
Let $K$ be a knot in $S^3$ and let $M_n\to S^3$ denote the branched $\Z/n\Z$-cover, 
that is, the Fox completion of the $\Z/n\Z$-cover $X_n\to X=S^3-K$. 
{In this setting, Livingston solved an analogue of Weber's problem in the study of knot concordance.}
\begin{theo}[{\cite[Theorem 1.2]{Livingston}}]
  The equality $|H_1(M_{p^n})_{\textup{tor}}|=1$ holds for every prime number $p$ and positive integer $n$ 
  if and only if 
  every non-trivial factor of the Alexander polynomial $\Delta_K(t)$ is 
  the $m$-th cyclotomic polynomial with $m$ being divisible by at least three distinct prime numbers. 
\end{theo}
He also showed that $M_{p^n}$ is a rational homology $3$-sphere for every knot $K$, prime $p$, and positive integer $n$.

On the other hand, the analogue of \cref{nconv} also holds.

\begin{theo}[{\cite[Theorem 3.1]{Ueki-Yoshizaki}}] \label{tconv}
  Let $(X_{p^n}\to X)_n$ be a $\Z_p$-cover of a compact 3-manifold $X$. 
  Then, the sizes of the torsion subgroups $H_1(X_{p^n})_{\textup{tor}}$, 
  those of the non-$p$ torsion subgroups $H_1(X_{p^n})_{\text{non-}p}$,
  and those of the $l$-torsion subgroups $H_1(X_{p^n})_{(l)}$ for each prime number $l$,  
  of the first homology groups converge in $\Z_p$. 
\end{theo}
The following lemma helps to prove the assertion in a parallel manner to \cref{nconv}. 
\begin{lem}[{\cite[Lemma 3.2]{Ueki-Yoshizaki}}] \label{surjsystem}
  Let $(X_{p^n}\to X)_n$ be a $\Z_p$-cover of a compact $3$-manifold.
  Then, $(H_1(X_{p^n})_{\textup{tor}})_n$ is a surjective system {for sufficiently large $n$.}
\end{lem} 

Recently Kionke \cite[Subsection 5.7]{Kionke} introduced new invariants, {\it the $p$-adic Betti number} and {\it the $p$-adic torsion}.
By Kionke's theorem \cite[Theorem 1.1 (ii)]{Kionke} and the Poincar\'{e} duality, the $p$-adic limit coincides with the $p$-adic torsion.

{Similarly to the case of function fields, we can compute the $p$-adic limits by using a well-known formula.}
{For a group $H$, we define $|H|$ by the size of $H$ if $H$ is a finite group and $0$ if otherwise.}
Let $\Delta_K(t)$ denote the Alexander polynomial of $K$ normalized by $\Delta_K(1)=1$.
Then we have
\begin{prop}[Fox--Weber's formula, cf.~\cite{CWeber}] \label{FoxWeber}
\[|H_1(M_n)|=|H_1(X_n)_{\textup{tor}}|=|\textup{Res}(t^n-1,\Delta_K(t))|\]
\end{prop}

\section{Examples}\label{Examples}
We present an explicit formula of the $p$-adic limits for polynomials in $\Z[t]$ {and apply it to} the cases of knots and function fields.

\subsection{Cyclic resultants of polynomials}\label{resultant}

Let $f(t)\in \Z[t]$.
We set
\[
  r_{p^n}(f)=\textup{Res}(t^{p^n}-1,f(t)).
\]
We note that if $f(t)$ is the Frobenius polynomial $F_k(t)$ of a function field $k$ or the Alexander polynomial $\Delta_K(t)$ of a knot $K$,
then we have $h(\F_{p^n}k)=|r_{p^n}(F_k(t))|$ or $|H_1(M_{p^n})|=|r_{p^n}(\Delta_K(t))|$ respectively, and $r_{p^n}$ converges to $\Z_p$.
In fact, {the convergence of $r_{p^n}(f)$ also hold for all polynomials $f(t)\in \Z[t]$.}
\begin{theo}[{\cite[Theorem 5.3]{Ueki-Yoshizaki}}]
  \label{thm.res.conv} 
  Let $0\neq f(t)\in \Z[t]$. 
  Then, the $p$-power-th cyclic resultants 
  $r_{p^n}(f)$ 
  converge in $\Z_p$. 
  The limit values are zero if and only if $p\mid f(1)$. 
  In any case, if $r_{p^n}(f) \neq 0$ for any $n$, then 
  the non-$p$-parts of $r_{p^n}(f)$ converge to a non-zero value in $\Z_p$. 
  For each prime number $l$, 
  similar assertions for the $l$-parts of $r_{p^n}(f)$ hold. 
\end{theo} 
The proof of the convergence is different from that of \cref{nconv}.
{A key step of the proof is to show $r_{p^n}(f)/r_{p^{n-1}}(f) = N_{\Q(\zeta_{p^n})/\Q}(f(\zeta_{p^n})) \equiv 1 \pmod {p^n}$ by using Artin reciprocity,
where $\zeta_p^n$ is a primitive $p^n$-th root of unity.}

In the following, we give an explicit formula for the $p$-adic limit of $r_{p^n}(f)$.
Let $\C_p$ denote the $p$-adic completion of an algebraic closure of the $p$-adic numbers $\Q_p$ and fix an embedding $\overline{\Q}\hookrightarrow \C_p$.
Let $\overline{\Z}_p$ denote the closure of $p$-adic integers $\Z_p$ in $\C_p$. 
Since extensions of $\F_p$ are $p$-prime-th cyclotomic extensions, 
an elementary $p$-adic number theory yields the following basic fact. 

\begin{lem}[{cf.~\cite[Lemma 2.10]{Ueki2}}] \label{lem.p-prime-th} 
If $\alpha\in \C_p$ satisfies $|\alpha|_p=1$, then there exists a unique $p$-prime-th root of unity $\zeta$ satisfying $|\alpha-\zeta|_p<1$. 
\end{lem} 

The following lemma is also elementary and classically known. 

\begin{lem} \label{lem.converge} 
Let $\alpha,\zeta \in \C_p$ 
with $|\alpha|_p=|\zeta|_p=1$. 

{\rm (1)} 
If $|\alpha-\zeta|_p<1$, then 
$\lim_{n\to \infty} \alpha^{p^n} -\zeta^{p^n} =0$ in $\C_p$. 

{\rm (2)} If $|\alpha-1|_p<1$, then 
$\lim_{n\to \infty} \frac{\alpha^{p^n} -1}{p^n}=\log \alpha$ in $\C_p$, 
where $\log$ denotes the $p$-adic logarithm defined by $\log (1+x)=\sum_{n=1}^\infty \frac{-(-x)^n}{n}$ on $\overline{\Z}_p$. 

{\rm (3)} If $|\alpha-1|_p<p^{-1/(p-1)}$, then $|\log \alpha|_p = |1-\alpha|_p$. 
\end{lem} 
By using \cref{lem.p-prime-th} and \ref{lem.converge}, we obtain the explicit formula for $p$-adic limits.
This theorem can be regarded as a refinement of \cite[Proposition 2]{Kisilevsky}.

\begin{theo}[{\cite[Theorem 5.7]{Ueki-Yoshizaki}}] \label{thm.res} 
  Let $0\neq f(t)\in \Z[t]$ and let $p^\mu$ denote the maximal $p$-power dividing $f(t)$. 
  Write $f(t)=a_0\prod_i (t-\alpha_i)$ in $\overline{\Q}[t]$ and note that $|p^{-\mu} a_0\prod_{|\alpha_i|_p>1} \alpha_i|_p=1$. 
  Let $\xi$ and $\zeta_i$ denote the unique $p$-prime-th roots of unity satisfying $|p^{-\mu}a_0\prod_{|\alpha_j|_p>1} \alpha_j -\xi|_p<1$ and $|\alpha_i-\zeta_i|_p<1$ for each $i$ with $|\alpha_i|_p=1$. 
  
  {\rm (1)} {\rm (i)} If $p\mid f(t)$, so that $\mu>0$, then $\lim_{n\to \infty} {\rm Res}(t^{p^n}-1,f(t)) =0$ holds in $\Z_p$.  
  
  {\rm (ii)} If $p\nmid f(t)$, so that $\mu=0$, then 
  \[\lim_{n\to \infty} {\rm Res}(t^{p^n}-1,f(t)) =
  (-1)^{p\, {\rm deg}f +\#\{i\, \mid\, |\alpha_i|_p<1\}} \xi \prod_{i;\,|\alpha_i|_p=1} (\zeta_i-1)\]
  holds in $\Z_p$, and the limit value is zero if and only if $\zeta_i=1$ for some $i$. 
  
  {\rm (2)} In any case, the non-$p$ part 
  ${\rm Res}(t^{p^n}-1,f(t))_{\text{non-}p}= {\rm Res}(t^{p^n}-1,f(t))\, | {\rm Res}(t^{p^n}-1,f(t))|_p$ converges to 
  \[(-1)^{p\, {\rm deg}f +\#\{i\, \mid\, |\alpha_i|_p<1\}} \xi\, \bigl( \prod_{\substack{i;\, |\alpha_i|_p=1,\\ \ |\alpha_i-1|_p=1}} (\zeta_i-1)\bigr) \, p^{-\nu} \prod_{\substack{i;\, |\alpha_i|_p=1,\\ \  |\alpha_i-1|_p<1}} \log \alpha_i\]
  in $\Z_p$, where $\log$ denotes the $p$-adic logarithm and 
  $\nu\in \Z$ is defined by 
  $p^{-\nu}=\prod_{\substack{i;\, 
  |\alpha_i-1|<1}} |\log \alpha_i|_p$. 
  If all $\alpha_i$'s with $|\alpha_i-1|_p<1$ are sufficiently close to 1, that is, 
  if they all satisfy $|\alpha_i-1|_p<p^{-1/(p-1)}$, then $p^\nu=|f(1)|_p^{-1}$ holds. 
\end{theo}
We note that if we put $\lambda=\#\{i\mid |\alpha_i-1|_p<1\}$, then these $\lambda, \mu, \nu$ are the Iwasawa invariants of $f(t)$, that is, $|{\rm Res}(t^{p^n}-1, f(t))|_p^{-1}=p^{\lambda n+\mu p^n +\nu}$ holds for any $n\gg 0$. 


\subsection{Function fields of elliptic curves}\label{ExamplesFF}
The Frobenius polynomial of a function field of genus $1$ is a monic quadratic polynomial.
More precisely, if $k$ is a function field of genus $1$ over the field $\F_q$ of $q$ elements, then we have
\[
  F_k(t)=t^2-a_1t+q,
\]
where we set $a_1=q+1-h(k)$.

{Let us examine a function field $k={\rm Frac}(\F_5[x,y]/(y^2-(x^3-1)))$ as an example.}
Since $h(k)=6$, we have $F_k(t)=t^2+5$.
We study $k_{p^n}=\F_{5^{p^n}}k$ for each prime number $p$.
By \cref{thm.res.conv}, if $p=2$ or $3$, then we see that $\lim_{n\to \infty}h(k_{p^n})=0$.
Indeed, we compute $h(k_{p^n})$ for $1\leq n \leq 3$ in the below by using PARI/GP.
\begin{table}[ht]
  \centering
  \begin{tabular}{c|c|c|c}
    $n$ & $1$  & $2$ & $3$ \\
    \hline
    $h(k_{2^n})$ & $2^2 3^2$  & $2^6 3^2$ & $2^8 3^2 13^2$ \\
    \hline
    $h(k_{3^n})$ & $2^1 3^2 7^1$  & $2^1 3^3 7^1 5167^1$ & $2^1 3^4 7^1 163^1 487^1 5167^1 16018507^1$\\
  \end{tabular}
\end{table}
We may observe that the $p$-part grows for each case.
By using \cref{thm.res} (2), we can compute the limit values of non-$p$ parts numerically.
In the case of $p=2$, since $|\sqrt{-5}-1|_2=1/2<1$ and $|-\sqrt{-5}-1|_2=1/2<1$, we have
\[
  \lim_{n\to\infty}|\Cl(k_{2^n})_{\textup{non-}2}|=\frac{\log (\sqrt{-5})\log(-\sqrt{-5})}{2^{\nu}}.
\]
Since the $2$-adic valuation of this limit is $0$,
and those of the numerator is $2$,  we see that $\nu=2$.
In the case of $p=3$, since $x^2+5\equiv (x+1)(x-1) \pmod 3$, the roots of unity satisfying \cref{lem.p-prime-th} is $1$ and $-1$.
Thus we obtain
\[
  \lim_{n\to\infty}|\Cl(k_{3^n})_{\textup{non-}3}|=-2\frac{\log (\sqrt{-5})}{3^{\nu}}.
\]
Here we choose a prime ideal $\mathfrak{p}=(3, \sqrt{-5}-1)$ over $3$ as an extension of the $3$-adic valuation to $\Q(\sqrt{-5})$.
By a similar calculation for the case $p=2$, we see that $\nu=1$.

Another interesting case is $p=5$.
By \cref{thm.res} (1) (ii), we have $\lim_{n\to \infty}h(k_{5^n})=1$.
This example answers the question asked in Subsection \ref{Inff}.
This is an example {that answer} the question we gave in \cref{Inff} when $p$-adic limits are $1$.

In fact, by our explicit formula, we can answer the question for function fields over prime field of genus $1$ except for $p=2$ and $3$.
Let $p$ be a prime number different from $2$ and $3$.
The limit value is $1$ if and only if
all roots of the Frobenius polynomial are not prime to $p$.
This is equivalent to $p=q$ and $p\mid a_1$.
The condition $p\mid a_1$ implies $a_1=0$ except for $p=2$ and $3$ by Hasse's bound $|a_1|<2\sqrt{q}$.
If $a_1=q+1-h(k)=0$, then the elliptic curve corresponding to $k$ is said to be supersingular.
By a basic complex multiplication (CM) theory, we can determine whether the elliptic curve over a prime field is supersingular (cf.~\cite[Theorem 12]{Lang}).
For a nonsingular algebraic curve $C_q$ over a finite field $\F_q$, let $k_{C_q}$ denote the corresponding function field.
In summary, we obtain the following.
\begin{prop}
    Let $E$ be an elliptic curve over $\Q$ with good reduction at $\F_q$ and suppose $q\neq 2,3$.
    Then we have $\lim_{n\to \infty} h(k_{E_{q^{q^n}}})=1$ in $\Z_q$ if and only if 
    $q$ is a supersinguler prime of $E$. 
    If in addition $E$ has 
    CM over $\Q(\sqrt{D})$ with $D<0$, 
    then $\lim_{n\to \infty}h(k_{E_{q^{q^n}}})=1$ holds in $\Z_q$ if and only if 
    Legendre's quadratic residue symbol
    satisfies $(\frac{D}{q})\neq1$. 
\end{prop}

\subsection{Knots}
Next, we observe knots.
First, we deal with torus knots.
Let $(a,b)$ be a coprime pair of integers. 
The Alexander polynomial 
\[\Delta_{K}=\dfrac{(1-t)(1-t^{ab})}{(1-t^a)(1-t^b)}=\prod_{\substack{m\mid ab\\ m\nmid a,\, m\nmid b}}\Phi_m(t)\]
of the $(a,b)$-torus knot $K=T_{a,b}$ is the product of cyclotomic polynomials.
\begin{figure}[H]
  \centering
  \includegraphics[width=55mm]{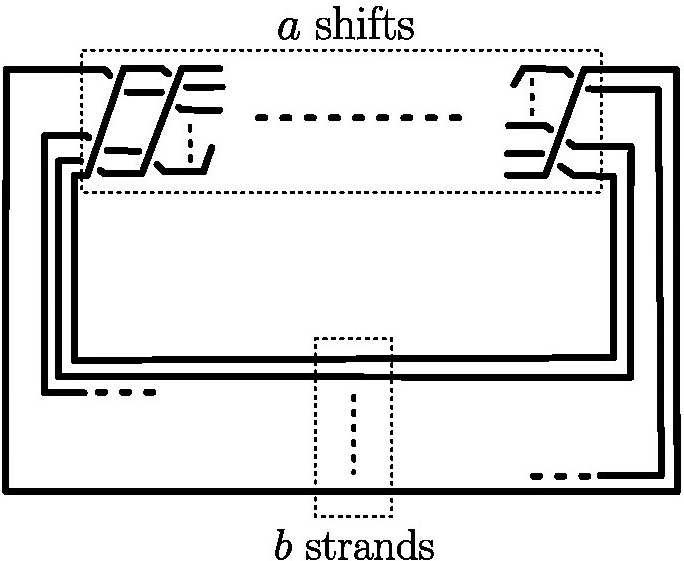}
  \caption*{$T_{a,b}$}
\end{figure}
For each $n\in \Z_{>0}$, let $\varphi(n)$ denote Euler's totient function. 
We invoke Apostol's result; 
\begin{lem}[{\cite[Theorem 4]{Apostol}}] \label{Lem.Apostol} 
Suppose that $m>n>1$ and $(m,n)>1$. Then,  
${\rm Res}(\Phi_m,\Phi_n)=p^{\varphi(n)}$ if $m/n$ is a power of a prime $p$, and 
${\rm Res}(\Phi_m,\Phi_n)=1$ if otherwise. 
\end{lem} 

\begin{prop}[{\cite[Proposition 6.3]{Ueki-Yoshizaki}}] \label{prop.torusknots}
Let $p$ be a prime number and let $(a,b)$ be a coprime pair of integers.
Assume that $p\nmid b$ and write $a=p^r a'$ with $r,a'\in \Z$, $p\nmid a'$. 
Let $(X_{p^n}\to X)_n$ denote the $\Z_p$-cover of the exterior of the torus knot $T_{a,b}$ in $S^3$. 
Then 
$|H_1(X_{p^n})_{\rm tor}|=b^{p^{{\rm min}\{n,r\}}-1}$ holds, and hence 
\[\lim_{n\to \infty}|H_1(X_{p^n})_{\rm tor}|=b^{p^r-1} \text{\ \ holds\ in\ }\Z_p.\] 
\end{prop}

For example,in the case of $(a,b)=(2,3)$, $K=T_{2,3}=J(2,2)=3_1$ (trefoil). Then we have $\Delta_K(t)=t^2-t+1=\Phi_6(t)$.
We have ${\rm Res}(t^{2^n}-1,\Delta_K(t))=3=3^{2-1}$, ${\rm Res}(t^{3^n}-1,\Delta_K(t))=4=2^{3-1}$, and ${\rm Res}(t^{p^n}-1,\Delta_K(t))=1=3^{p^0-1}$ for $p\neq 2,3$ for all $n\in \Z_{>0}$. 

Secondly, we deal with twist knots.
For each $m\in \Z$, the Alexander polynomial of the twist knot $K=J(2,2m)$ is given by 
$\Delta_{J(2,2m)}(t)=mt^2+(1-2m)t+m$. 
The convention is due to \cite{HosteShanahan2004JKTR}, so that we have $J(2,2)=3_1$ and $J(2,-2)=4_1$ for instance. 
\begin{figure}[H]
  \centering
  \includegraphics[width=45mm]{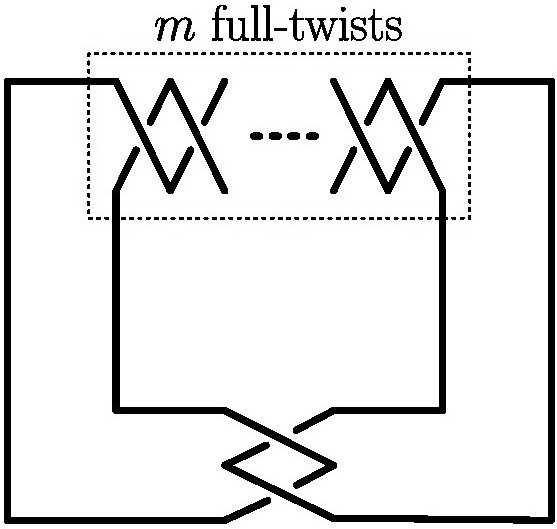}
  \caption*{$J(2,2m)$}
\end{figure}
Let $K=J(2,6)=7_2$. Then we have $\Delta_K(t)=3t^2-5t+3$ and 
\begin{center}
\begin{tabular}{c|c|c|c|c|c|c|c}
$p$&2&3&5&7&$\cdots$ \ \\ \hline 
$\lim_{n\to \infty}|H_1(X_{p^n})_{\rm tor}|$&$3$&$1$&$-2\sqrt{-1}$&$(2+\sqrt{2})\frac{1-\sqrt{-3}}{2}$&$\cdots$ \ 
\end{tabular}.
\end{center} 
Indeed, we compute the $7$-adic expansions of $|H_1(X_{7^6})_{\rm tor}|$ and $(2+\sqrt{2})(1-\sqrt{-3})/2$ by using PARI/GP as follows;
\begin{align*}
  |H_1(X_{7^6})_{\rm tor}|&=1 + 4*7 + 7^2 + 4*7^3 + 6*7^5 + 6*7^6 + 2*7^7 + 2*7^8... \\
  (2+\sqrt{2})(1-\sqrt{-3})/2&=1 + 4*7 + 7^2 + 4*7^3 + 6*7^5 + 6*7^6 + 6*7^7 + 7^8...
\end{align*}
This example shows that these two elements are close in $\Z_7$.

For twist knots, we can find all cases with the $p$-adic limits in $\Z$.
\begin{prop}[{\cite[Proposition 6.8, 6.10]{Ueki-Yoshizaki}}]
  Let us consider the $\Z_p$-cover $(X_{p^n}\to X)$ of the exterior of the twist knot $K=J(2,2m)$.
  \begin{itemize}
    \item[\textup{(1)}] If $p\mid m$, then we have
    \[\lim_{n\to \infty}|H_1(X_{p^n})_{\rm tor}|=
    \begin{cases}
    -{\rm sgn}(4m-1)&\text{if}\ p=2,\\
    1&\text{if}\ p\neq 2
    \end{cases} \text{\ \ in\ }\Z_p.\]
    \item[\textup{(2)}] If $p\nmid m$, then we have $\lim=\lim_{n\to \infty}|H_1(X_{p^n})_{\rm tor}| \in \Z$ if and only if one of the following holds.
    \begin{itemize}
      \item[$\bullet$] $p=2${\rm ;} $\lim={\rm sgn}(4m-1)\cdot 3=\pm3$.
      \item[$\bullet$] $p=3${\rm ;} $3\mid m-1$, $\lim=4$ or $3\mid m+1$, $\lim=-2$.
      \item[$\bullet$] $p=5${\rm ;} $5\mid m+1$, $\lim=-4$.
      \item[$\bullet$] $p\neq 2$, $3${\rm ;} $p\mid m-1$, $\lim=1$.
    \end{itemize}
  \end{itemize}
\end{prop}
By this result, we may conclude the following, answering a weak analogue of Weber's class number problem.
\begin{cor}[{\cite[Corollary 6.11]{Ueki-Yoshizaki}}]
  We have $\lim_{n\to \infty} |H_1(X_{p^n})|=1$ if and only if the following holds.
  \begin{itemize}
    \item $p=2$, $m$ is even, and $4m-1<0$.
    \item $p\neq 2$ and $p\mid m$.
    \item $p\neq 2$, $3$ and $p\mid (m-1)$.
  \end{itemize}
\end{cor}

\subsection{Remarks}
For number fields, there are few numerical results on the $p$-adic limits of class numbers in $\Z_p$-extensions.
Nevertheless, Ozaki \cite{Ozaki} studies the relation between several invariants of number fields and the $p$-adic limits.

For function fields, we have not yet been able to compute the $p$-adic limits of geometric $\Z_p$-extensions.

For knots, we studied torus knots and twist knots as concrete examples.
The cases of links and multivariable polynomials will be discussed in \cite{Tateno-Ueki}.

Our numerical study of the $p$-adic limits of class numbers for knots gives explicit examples of Kionke's $p$-adic torsions \cite{Kionke}.
Kionke also discusses whether the $p$-adic Betti number belongs to $\Z$, as the $p$-adic analogue of Atiyah conjecture.
It seems to be interesting to ask when the $p$-adic torsion belongs to $\Z$.

By changing the base field (resp. space) $k$ to $k_e$ (resp. $X$ to $X_e$) for some integer $e$ with $p \nmid e$,
we may pursue systematic numerical studies of the $p$-adic limits and Iwasawa $\nu$-invariants as well.
The article further concerns with the following questions:
(1) When does the $p$-adic limit belong to $\Z$ ?
(2) Can $\nu$ be arbitrarily large, with the base $p$-torsion being small? 
For twist knots and elliptic curves over $\F_l$, we have complete lists for (1) \cite[Proposition 6.10, Proposition 7.10]{Ueki-Yoshizaki} and answers to (2) in affirmative with precise classifications \cite[Example 6.14, 6.15, 7.15, 7.16, Proposition 6.16]{Ueki-Yoshizaki}.


As we mentioned in the first section, the unit group is an essential object in number fields.
However, in topology, the studies on analogues of the unit group have not progressed as much as in number theory.
Among other things, the analytic class number formula is an important connection between the class groups and the unit groups.
Therefore, finding an analogue of the analytic class number formula in topology would be a very important problem.

\section*{Acknowledgement}
The author would like to express his sincere gratitude to
Tomokazu Kashio, Manabu Ozaki, Makoto Sakuma, Sohei Tateno, and Jun Ueki
for useful information and fruitful conversations.
The author has been supported by JSPS KAKENHI Grant Number 22J10004.

\bibliographystyle{alpha}
\bibliography{references}

\begin{thebibliography}{FKM14}

\bibitem[AP04]{Aubry-Perret}
Yves Aubry and Marc Perret.
\newblock {On the characteristic polynomials of the Frobenius endomorphism for
  projective curves over finite fields}.
\newblock {\em {Finite Fields and Their Applications}}, 10(3):412--431, 2004.

\bibitem[Apo70]{Apostol}
Tom~M. Apostol.
\newblock {Resultants of Cyclotomic Polynomials}.
\newblock {\em Proc. Amer. Math. Soc.}, 24(3):457--462, 1970.

\bibitem[Bau69]{Bauer}
Helmut Bauer.
\newblock {Numerische Bestimmung von Klassenzahlen reeller zyklischer
  Zahlk\"{o}rper}.
\newblock {\em Journal of Number Theory}, 1(2):161--162, 1969.

\bibitem[FK09]{Fukuda-Komatsu1}
Takashi Fukuda and Keiichi Komatsu.
\newblock {Weber's Class Number Problem in the Cyclotomic $\Z_2$-Extension of
  $\mathbb{Q}$}.
\newblock {\em Experimental Mathematics}, 18(2):213--222, 2009.

\bibitem[FK11]{Fukuda-Komatsu3}
Takashi Fukuda and Keiichi Komatsu.
\newblock Weber's class number problem in the cyclotomic $\z_2$-extension of
  $\mathbb{Q}$, iii.
\newblock {\em International Journal of Number Theory}, 07(06):1627--1635,
  2011.

\bibitem[FKM14]{FKM}
Takashi Fukuda, Keiichi Komatsu, and Takayuki Morisawa.
\newblock Weber's class number one problem.
\newblock In Thanasis Bouganis and Otmar Venjakob, editors, {\em Iwasawa Theory
  2012}, pages 221--226, Berlin, Heidelberg, 2014. Springer Berlin Heidelberg.

\bibitem[Fuk19]{Fukuda-Kaisetsu}
Takashi Fukuda.
\newblock {\em Juuten kaisetsu Iwasawa riron: riron kara keisan made (in
  Japanese)}.
\newblock Number 145 in SGC Library. SAIENSU - SHA, 2019.

\bibitem[Hor05]{Horie05}
Kuniaki Horie.
\newblock {The ideal class group of the basic $\boldsymbol{Z}_p$-extension over
  an imaginary quadratic field}.
\newblock {\em Tohoku Mathematical Journal}, 57, 09 2005.

\bibitem[Hor07]{Horie07}
Kuniaki Horie.
\newblock {Certain primary components of the ideal class group of the
  $\boldsymbol{Z}_p$-extension over the rationals}.
\newblock {\em Tohoku Mathematical Journal}, 59(2):259 -- 291, 2007.

\bibitem[HS04]{HosteShanahan2004JKTR}
Jim Hoste and Patrick~D. Shanahan.
\newblock {A} formula for the {A}-polynomial of twist knots.
\newblock {\em Journal of Knot Theory and Its Ramifications}, 13(02):193--209,
  2004.

\bibitem[Iwa63]{Iwasawa}
Kenkichi Iwasawa.
\newblock {Daisutai to kansutai tono aru ruiji ni tsuite (in Japanese)}.
\newblock {\em Sugaku}, 15(2):65--67, 1963.

\bibitem[Kio20]{Kionke}
Steffen Kionke.
\newblock On $p$-adic limits of topological invariants.
\newblock {\em Journal of the London Mathematical Society}, 102(2):498--534,
  apr 2020.

\bibitem[Kis97]{Kisilevsky}
Hershy Kisilevsky.
\newblock {A Generalization of a result of Sinnott}.
\newblock {\em Pacific Journal of Mathematics}, 181(3):225--229, 1997.

\bibitem[KY21]{Kashio-Yoshizaki}
Tomokazu Kashio and Hyuga Yoshizaki.
\newblock {Minimal relative units of the cyclotomic $\mathbb{Z}_2$-extension},
  2021.
\newblock preprint. arXiv:2107.08587.

\bibitem[Lan87]{Lang}
Serge Lang.
\newblock {\em {Elliptic Functions}}.
\newblock Graduate texts in mathematics. Springer, 1987.

\bibitem[Liv02]{Livingston}
Charles Livingston.
\newblock {Seifert forms and concordance}.
\newblock {\em Geometry \& Topology}, 6(1):403--408, 2002.

\bibitem[Mas78]{Masley}
John~Myron Masley.
\newblock {Class numbers of real cyclic number fields with small conductor}.
\newblock {\em Compositio Mathematica}, 37(3):297--319, 1978.

\bibitem[Maz64]{Mazur}
Barry Mazur.
\newblock {Remarks on the Alexander polynomial}.
\newblock http://www.math.harvard.edu/~mazur/10 papers/alexander
  polynomial.pdf, 1963--64.

\bibitem[Mil14]{Miller}
John~C. Miller.
\newblock {Class numbers of totally real fields and applications to the Weber
  class number problem}.
\newblock {\em Acta Arith.}, 164:381--397, 2014.

\bibitem[MO16]{Morisawa-Okazaki}
Takayuki Morisawa and Ryotaro Okazaki.
\newblock {Height and Weber’s Class Number Problem}.
\newblock {\em Journal de Th{\'e}orie des Nombres de Bordeaux}, Tome
  28(3):811--828, 2016.

\bibitem[MO20]{MO}
Takayuki Morisawa and Ryotaro Okazaki.
\newblock {Filtrations of units of Vi{\`e}te field}.
\newblock {\em International Journal of Number Theory}, 16(05):1067--1079,
  2020.

\bibitem[MS15]{Mercuri-Stirpe}
Pietro Mercuri and Claudio Stirpe.
\newblock Classification of algebraic function fields with class number one.
\newblock {\em Journal of Number Theory}, 154:365--374, 2015.

\bibitem[Oza22]{Ozaki}
Manabu Ozaki.
\newblock {On the $p$-adic limit of class numbers along a pro-$p$-extension}.
\newblock preprint, 2022.

\bibitem[Ros02]{Rosen}
Michael Rosen.
\newblock {\em Number Theory in Function Fields}.
\newblock Graduate Texts in Mathematics. Springer New York, 2002.

\bibitem[SS15]{Shen-Shi}
Qibin Shen and Shuhui Shi.
\newblock {Function fields of class number one}.
\newblock {\em Journal of Number Theory}, 154:375--379, 2015.

\bibitem[TU22]{Tateno-Ueki}
Sohei Tateno and Jun Ueki.
\newblock {The Iwasawa invariants of $\mathbb{Z}_p^d$ covers of links}, 2022.
\newblock preprint. arXiv:.

\bibitem[Uek20]{Ueki2}
Jun Ueki.
\newblock {$p$-adic Mahler measure and $\mathbb{Z}$-covers of links}.
\newblock {\em Ergodic Theory and Dynamical Systems}, 40(1), 2020.

\bibitem[UY22]{Ueki-Yoshizaki}
Jun Ueki and Hyuga Yoshizaki.
\newblock The $p$-adic limits of class numbers in $\mathbb{Z}_p$-towers, 2022.
\newblock preprint. arXiv:2210.06182.

\bibitem[vdL82]{Linden}
F.~J. van~der Linden.
\newblock {Class Number Computations of Real Abelian Number Fields}.
\newblock {\em Mathematics of Computation}, 39(160):693--707, 1982.

\bibitem[Was97]{Washington}
Lawrence~C. Washington.
\newblock {\em Introduction to Cyclotomic Fields}.
\newblock Springer New York, NY, 2 edition, 1997.

\bibitem[Web00]{HWeber}
Heinrich Weber.
\newblock {Theorie der Abel'schen Zahlk\"{o}rper}.
\newblock {\em Acta Mathematica}, 8:193 -- 263, 1900.

\bibitem[Web80]{CWeber}
Claude Weber.
\newblock {Sur une formule de R. H. Fox concernant l'homologie des
  rev\^etements cycliques}.
\newblock {\em Enseign. Math. (2)}, 25:261, 1980.

\bibitem[Wei48]{Weil}
Andr\'{e} Weil.
\newblock {\em {Sur les courbes alg{\'e}briques et les vari{\'e}t{\'e}s qui
  s'en d{\'e}duisent}}.
\newblock Actualit{\'e}s scientifiques et industrielles. Hermann, 1948.

\bibitem[Yos20]{Yoshizaki}
Hyuga Yoshizaki.
\newblock {Generalized Pell's equations and Weber's class number problem},
  2020.
\newblock preprint. arXiv:2010.06399.

\end{thebibliography}

\end{document}